\documentclass[12pt, reqno]{amsart}
\usepackage{latexsym,epsfig,graphpap,graphics,amsmath,amssymb, amsthm,  amsfonts, color}
\usepackage[bookmarksnumbered, plainpages, backref]{hyperref}
\numberwithin{equation}{section}

\begin{document}
\setcounter{page}{1}
\begin{center}
\bf{WARPED PRODUCT SUBMANIFOLDS OF KENMOTSU MANIFOLDS}
\end{center} 
\begin{center} Sachin Kumar Srivastava\end{center}
\vspace{.3cm}\noindent {\bf Abstract:}~ In the present paper, we study warped product   semi-slant submanifolds of Kenmotsu manifolds.We have obtained results on the existence of warped product semi-slant submanifolds of Kenmotsu manifolds in term of the canonical structure.\\
\vspace{.3cm}\noindent{\bf{2000 Mathematics Subject Classification :}}~{53C40, 53C42, 53B25.}\\
\noindent{\bf{Keywords :}}~{Warped product; doubly warped product; slant submanifold; semi-slant submanifold; Kenmotsu manifolds; canonical structure.}

\section{Introduction}

\indent The differential geometry of slant submanifolds has shown an increasing development since B. - Y. Chen defined slant immersion in complex geometry as a natural generalization of both holomorphic immersions and totally real immersions (see \cite{BYC}).In (\cite{AL183}), A. Lotta has introduced the notion of slant immersion of a Riemannian manifold into an almost contact metric manifold. In (\cite{RLB1}), the notion of warped product manifolds was introduced by Bishop and O'Neill in 1969 and it was studied by many mathematicians and physicists.\newline
\indent  The notion of semi-slant submanifolds of almost Hermitian manifolds was introduced by N. Papaghuic (see \cite{NP55}). In fact, semi-slant submanifolds in almost Hermitian manifolds are defined on the line of CR-submanifolds.  In the setting of almost contact metric manifolds, semi-slant submanifolds are defined and investigated by J.L. Cabrerizo et. al (see \cite{JLC125}).

\indent The study of warped product semi-slant submanifolds of Kaehler manifolds was introduced by B. Sahin (see \cite{BS195}). Later, K.A. Khan et.al studied   warped product   semi-slant submanifolds in cosymplectic manifolds and showed that there exist no proper warped product semi-slant submanifolds in the forms $N_{T} \times _{f} N_{\theta }$ and reversing the two factors in cosymplectic manifolds (see \cite{KA55}).

\indent Recently, M. Atceken proved that the warped product submanifolds of the types $M=N_{T} \times _{f} N_{\theta }$ and $M=N_{\theta } \times _{f} N_{\bot } $of a Kenmotsu manifold $\tilde{M}$do not exist where the manifolds $N_{\theta } $   and $N_{T} $  (resp.$N_{\bot } $) are proper slant and $\varphi $-invariant (resp.anti-invariant) submanifolds of a Kenmotsu manifold  $\tilde{M}$,  respectively (see \cite{MAE1}).  In this paper, we have obtained some results for the existence of warped product semi-slant submanifolds of $\beta-$ Kenmotsu manifolds. 

\section{Preliminaries}

\indent Let $\tilde{M}$ be a (2n + 1) -- dimensional $C^{\infty }$ manifold endowed with the almost contact metric structure$\left(\phi ,\xi ,\eta ,g\right)$, where $\phi$ is a tensor field of type (1, 1), $\xi$ is a vector field, $\eta$ is a 1 -- form and g is a Riemannian metric on  $\tilde{M}$, all these tensor fields satisfying (\cite{DEB})

\begin{equation} \label{GrindEQ__2_1_}\varphi ^{2} =I-\eta \otimes \xi ,\eta \left(\xi \right)=1,\eta o\varphi =0, \varphi \xi =0 \end{equation}

 \begin{equation} \label{GrindEQ__2_2_}\left(X,Y\right)=g\left(\varphi X,\varphi Y\right)+\eta \left(X\right)\, \eta \left(Y\right)\end{equation}

\begin{equation}\label{GrindEQ__2_3_} g\left(X,\varphi Y\right)=-g\left(\varphi X,Y\right), g\left(X,\xi \right)=\eta \, \left(X\right)
\end{equation}
\noindent for any $X,Y\in T\tilde{M}$. Here $T\tilde{M}$is the standard notation for the tangent bundle of $\tilde{M}$. The two forms $\Phi $ denote the fundamental two forms and is given by $g\left(X,\phi Y\right)=\Phi \left(X,Y\right)$. The manifold is said to be contact if $\Phi =d\eta $. If $\xi $ is a killing vector field with respect to g, the contact metric structure is called K - contact structure. It is known that a contact metric manifold is K -- contact if and only if $\tilde{\nabla }_{X} \xi =-\phi X$, where $\tilde{\nabla }$ denotes the Levi -- Civita connection on $\tilde{M}$.\\
\indent The almost contact structure $\tilde{M}$ is said to be normal if $[\phi ,\phi ]+2d\eta \otimes \xi =0$, where $[\phi ,\phi ]$ is the Nijenhuis tensor of $\phi $. A Sasakian manifold is a normal contact metric manifold. Every Sasakian manifold is K -- contact. A three dimensional K -- contact manifold is Sasakian. An almost contact metric manifold is Sasakian if and only if 
\begin{equation} \label{GrindEQ__2_4_}
\left(\tilde{\nabla }_{X} \varphi \right)Y=g\left(X,Y\right)\xi -\eta \left(Y\right)X;  X,Y\in T\tilde{M}
\end{equation}
\noindent Moreover, on a Sasakian manifold 
\begin{equation} \label{GrindEQ__2_5_}\tilde{\nabla }_{X} \xi =-\varphi X.
\end{equation}
\noindent For any $X\in T\tilde{M}$and $\xi $ is the structure vector field.
\noindent An almost contact metric structure $\left(\phi ,\xi ,\eta ,g\right)$ on $\tilde{M}$ is called $\beta- $ Kenmotsu manifold if 

\begin{equation} \label{GrindEQ__2_6_}
\left(\tilde{\nabla }_{X} \varphi \right)Y=\beta \left\{g\left(\varphi X,Y\right)\xi -\eta \left(Y\right)\varphi X\right\}
\end{equation}
\noindent where $\beta $ is a smooth function on $\tilde{M}$ and $\tilde{\nabla }$means the covariant differentiation with respect to g (see \cite{JO187}).

\noindent If we put $\beta =1$ in \eqref{GrindEQ__2_6_} then $\tilde{M}$is a Kenmotsu manifold. Let M be a submanifold immersed in a (2n + 1) -- dimensional contact metric manifold $\tilde{M}$; we denote by the same symbol g the induced metric on M. TM is the tangent bundle of the manifold M and $T^{\bot } M$ is the set of vector fields normal to M. Then the Gauss and Weingarten formula is given by

\begin{equation} \label{GrindEQ__2_7_}\tilde{\nabla }_{X} Y=\nabla _{X} Y+h(X,Y) \end{equation}

\begin{equation} \label{GrindEQ__2_8_}\tilde{\nabla }_{X} N=-A_{N} X+\nabla _{X}^{\bot } Y\end{equation}

\noindent For any $X,Y\in TM$ and $N\in T^{\bot } M$, where $\nabla ^{\bot } $is the connection in the normal bundle. The second fundamental form h and the shape operator $A_{N} $ are related by

\begin{equation} \label{GrindEQ__2_9_} g\left(A_{N} X,Y\right)=g\left(h(X,Y),N\right).\end{equation}

\noindent For any $X,Y\in TM$ and $N\in T^{\bot } M$, we write 

\begin{equation} \label{GrindEQ__2_10_}\varphi X=TX+NX,\, \left(TX\in TM,NX\in T^{\bot } M\right),\end{equation}
\begin{equation} \label{GrindEQ__2_11_}\varphi N=tN+nN,\, \left(tN\in TM,nN\in T^{\bot } M\right)\,.\end{equation}

\noindent The submanifold M is invariant if N is identically zero on the other hand; M is anti -- invariant if T is identically zero. From \eqref{GrindEQ__2_3_} and \eqref{GrindEQ__2_10_}, we have 

\begin{equation} \label{GrindEQ__2_12_} g\left(X,TY\right)=-g\left(TX,Y\right)\end{equation}
\noindent For any $X,Y\in TM$.
\noindent The covariant derivatives of the tensor fields T and N are defined as
\begin{equation}\label{GrindEQ__2_13_} \left(\nabla _{X} T\right)Y=\nabla _{X} TY-T\nabla _{X}Y\end{equation}
\begin{equation}\label{GrindEQ__2_14_} (\tilde{\nabla }_{X} N)Y=\nabla _{X}^{\bot } NY-N\nabla _{X}Y \end{equation}
\noindent For all $X,Y\in TM$.
\noindent The canonical structure T and N on a submanifold M are said to be \textit{parallel}if  $\nabla T=0$ and $\tilde{\nabla }N=0$, respectively.\\
\indent We shall always consider $\xi $ to be tangent to M. The submanifold M is said to be \textit{invariant} if N is identically zero, that is, $\phi X\in TM$ for any $X\in TM$.On the other hand M is said to be \textit{anti -- invariant} if T is identically zero, that is, $\phi X\in T^{\bot } M$ for any$X\in TM$. For each non -- zero X tangent to M, such that X is not proportional to $\xi $, we denote by $\theta \left(X\right)$, the angle between $\phi $ and TX.
\noindent M is said to be \textit{slant} \cite{JLC125} if the angle $\theta \left(X\right)$is constant for all $X\in TM-\left\{\xi \right\}$. The angle $\theta$ is called \textit{slant angle} or \textit{Wiritinger angle}. Obviously if $\theta =0$, M is invariant and if $\theta =\pi /2$, M is anti -- invariant submanifold. If the slant angle of M is different from 0 and $\pi /2$ then it is called \textit{proper slant}.\\
\noindent A characterization of slant submanifolds is given by the following:\\
\noindent \textbf{Theorem 2.1. (\cite{JLC125})} \textit{Let M be a submanifold of an almost contact metric manifold M, such that $\xi \in TM$. Then M is slant if and only if there exists a constant $\lambda \in [0,1]$ such that}
\begin{equation} \label{GrindEQ__2_15_}T^{2} =\lambda \left(-I+\eta \otimes \xi \right)\end{equation}
\noindent \textit{Furthermore, if $\theta$ is slant angle then $\lambda =\cos ^{2} \theta $.}\\
\noindent The following relations are straight forward consequences of equation\eqref{GrindEQ__2_15_}:
\begin{equation} \label{GrindEQ__2_16_} g\left(TX,TY\right)=\cos ^{2} \theta \, \left[g\left(X,Y\right)-\eta \left(X\right)\, \eta \left(Y\right)\right] \end{equation}
\begin{equation} \label{GrindEQ__2_17_} g\left(NX,NY\right)=\sin ^{2} \theta \, \left[g\left(X,Y\right)-\eta \left(X\right)\, \eta \left(Y\right)\right] \end{equation}
\noindent for any $X,Y\in TM$.
\noindent We say M is a semi -- slant submanifold of  $\tilde{M}$ if there exist an orthogonal direct decomposition of TM as $TM=D_{\, 1} \oplus D_{\, 2} \oplus \left\{\xi \right\}$ , where $D_{\, 1} $is invariant distribution i.e., $\phi \left(D_{\, 1} \right)=D_{\, 1}$ and $D_{\, 2}$ is slant with slant angle $\theta \ne 0.$ The orthogonal complement of  $ND_{\, 2} $is the normal bundle $T^{\bot } M$, is an invariant subbundle of  $T^{\bot }M$ and is denoted by $\mu $. Thus, we have $T^{\bot } M=ND_{2} \oplus \mu$.
\noindent Similarly, we say that M is anti -- slant submanifold of $\tilde{M}$ if $D_{1} $is an anti -- invariant distribution of M i.e., $\phi D_{1} \subseteq T^{\bot } M$and $D_{2}$ is slant with slant angle $\theta \ne 0$.
\section{Warped and doubly product manifolds}
\noindent Let $\left(N_{1} ,g_{1} \right)$and $\left(N_{2} ,g_{2} \right)$ be two Riemannian manifolds and {\it f}, a positive differentiable function on $N_{1} $. The warped product of $N_{1} $and $N_{2} $ is the Riemannian manifold $N_{1} \times {}_{f} N_{2} =\left(N_{1} \times N_{2} ,g\right)$, where
\begin{equation} \label{GrindEQ__3_1_} g=g_{1} +f^{2} g_{2}.
\end{equation}
\noindent A warped product manifold $N_{1} \times {}_{f} N_{2} $ is said to be trivial if the warping function {\it f} is constant. We recall the following general formula on warped product (see \cite{RLB1}).

\begin{equation} \label{GrindEQ__3_2_}\nabla _{X} V=\nabla _{V} X=\left(X\ln f\right)V,
\end{equation}

\noindent where X is tangent to $N_{1} $and V is tangent to $N_{2}$.

\noindent Let $M=N_{1} \times {}_{f} N_{2} $ be a warped product manifold, this means that $N_{1}$ is totally geodesic and $N_{2} $ is totally umbilical submanifold of M, respectively.
\noindent
Doubly warped product manifolds were introduced as a generalization of warped product manifolds by B. Unal \cite{BU253}. A doubly warped product manifold of $N_{1} $and$N_{2} $, denoted as ${}_{f_{2} } N_{1} \times {}_{f_{1} } N_{2} $ is endowed with a metric g defined as 
\begin{equation} \label{GrindEQ__3_3_} g=f_{2} ^{2} g_{1} +f_{1} ^{2} g_{2} 
\end{equation}
\noindent 
Where $f_{1} $ and $f_{2} $ are positive differentiable functions on $N_{1} $ and $N_{2} $ respectively. In this case formula \eqref{GrindEQ__3_2_} is generalized as
\begin{equation} \label{GrindEQ__3_4_}\nabla _{X} V=\left(X\ln f_{1} \right)V+\left(V\ln f_{2} \right)X
\end{equation}
\noindent 
for each, $X\in TN_{1} $ and $V\in TN_{2}$ ({see \cite{MIM121}). 
\noindent
If neither $f_{1} $ nor $f_{2} $ is constant we have a non trivial doubly warped product $ M=N_{1} \times {}_{f} N_{2}$. Obviously in this case both $N_{1}$ and $N_{2}$ are totally umbilical submanifolds of M.

\noindent For any $X\in TN_{1}$ and $Z\in TN_{2}$ then by \eqref{GrindEQ__3_4_}, we have 
\[\nabla _{X} Z=\nabla _{Z} X=\left(X\ln f_{1} \right)Z+\left(Z\ln f_{2} \right)X.\] 
If $\xi \in TN_{1}$ then above equation gives

\begin{equation} \label{GrindEQ__3_5_}\nabla _{\xi } Z=\nabla _{Z} \xi =\left(\xi \ln f_{1} \right)Z+\left(Z\ln f_{2} \right)\xi.
\end{equation}

\noindent On the other hand, using equation \eqref{GrindEQ__2_6_} and the fact that $\xi $ is tangent to $N_{1} $, we have
\[\tilde{\nabla }_{Z} \xi =\beta Z.\] 
Using equation {\eqref{GrindEQ__2_7_}, we have 
\[\nabla _{Z} \xi +h\left(Z,\xi \right)=\beta Z.\] 
\noindent Using equation{\eqref{GrindEQ__3_5_} and then comparing the tangential and normal component we obtain 
\begin{equation} \label{GrindEQ__3_6_} \left(\xi \ln f_{1} \right)Z+\left(Z\ln f_{2} \right)\xi =\beta Z,\end{equation}
\begin{equation} \label{GrindEQ__3_7_}h (Z,\xi ) = 0. \end{equation}
\noindent Taking product with Z in equation \eqref{GrindEQ__3_6_} and using the fact that $\xi $, Z and TZ are orthogonal vector fields then$\left(\xi \ln f_{1} \right)=1,\, \, \left(Z\ln f_{2} \right)=0$.This shows that $f_{2} $ is constant.
\noindent Similarly, if the structure vector field $\xi$ is tangent to $N_{2} $ and for any $X\in TN_{1} $ we obtain$\left(\xi \ln f_{2} \right)=1,\, \, \left(X\ln f_{1} \right)=0$. Showing that $f_{1} $ is constant.
\noindent 
This leads to the following theorem:\\
\noindent \textbf{Theorem 3.1.} \textit{There do not exist proper doubly warped product submanifolds $M={}_{f_{2} } N_{1} \times {}_{f_{1} } N_{2} $ of a Kenmotsu manifold $\tilde{M}$ where $N_{1} $ and $N_{2} $ are any Riemannian submanifolds of $\tilde{M}$}.\\
\noindent The following corollary is an immediate consequence of the above theorem:\\
\noindent \textbf{Corollary 3.1.} \textit{There do not exist warped product submanifolds $M=N_{1} \times {}_{f} N_{2} $ of a Kenmotsu manifold $\tilde{M}$ such that $\xi \in TN_{2} $ where $N_{1} $ and $N_{2} $ are any Riemannian submanifolds of $\tilde{M}$}.\\
\noindent Thus, the only remaining case to study of warped product submanifolds $M=N_{1} \times {}_{f} N_{2} $ of a Kenmotsu manifold $\tilde{M}$ such that $\xi \in TN_{1} $. For any $X\in TN_{1} $ and $Z\in TN_{2}$, we have
\[\left(\tilde{\nabla }_{X} \phi \right)=\tilde{\nabla }_{X} \phi Z-\phi \tilde{\nabla }_{X} Z.\] 
\noindent Using equation {\eqref{GrindEQ__2_6_}and the fact that $\xi \in TN_{1}$, left hand side of the above equation is zero by orthogonality of two distributions, then $\phi \tilde{\nabla }_{X} Z=\tilde{\nabla }_{X} \phi Z$.
\noindent By equations \eqref{GrindEQ__2_7_}, \eqref{GrindEQ__2_8_}, \eqref{GrindEQ__2_10_} and \eqref{GrindEQ__2_11_}, we obtain
\[\nabla _{X} TZ+h\left(X,TZ\right)-A_{NZ} X+\nabla _{X}^{\bot } NZ=T\nabla _{X} Z+N\nabla _{X} Z+th\left(X,Z\right)+nh\left(X,Z\right).\] 
\noindent Equating the tangential and normal components and using \eqref{GrindEQ__3_2_}, we get
\begin{equation} \label{GrindEQ__3_8_} A_{NZ} X=-th\left(X,Z\right)
\end{equation} 
\noindent and

\begin{equation} \label{GrindEQ__3_9_}\nabla _{X}^{\bot } NZ=\left(X\ln f\right)NZ+nh\left(X,Z\right)-h\left(X,TZ\right).
\end{equation}
\noindent
This leads to the following theorem:\\
\noindent \textbf{Theorem 3.2.} \textit{Let $M=N_{1} \times {}_{f} N_{2} $ be warped product submanifolds of a Kenmotsu manifold $\tilde{M}$ such that $\xi \in TN_{1} $ where $N_{1} $ and $N_{2} $ are any Riemannian submanifolds of $\tilde{M}$.Then}

\begin{equation} \label{GrindEQ__3_10_}\xi \ln f=1,
\end{equation}
\begin{equation}\label{GrindEQ__3_11_} A_{NZ}X=-Bh\left(X,Z\right),
\end{equation} 
\begin{equation} \label{GrindEQ__3_12_} g\left(h\left(X,Y\right),NZ\right)=g\left(h\left(X,Z\right),NY\right),
\end{equation}
\begin{equation} \label{GrindEQ__3_13_} g\left(h\left(X,Z\right),NW\right)=g\left(h\left(X,W\right),NZ\right)
\end{equation}
\noindent
\textit{for any $X,Y\in TN_{1} $ and $Z,W\in TN_{2} $}.
\section{Warped product semi-slant submanifolds}

\noindent 

\noindent We have seen that the warped products   of the type $N_{1} \times {}_{f} N_{2} $ of Kenmotsu manifolds do not exist if \textit{$\xi \in TN_{2} $}.Thus, in this section we study warped product semi -- slant submanifolds $N_{1} \times {}_{f} N_{2}$ of Kenmotsu manifolds only when \textit{$\xi \in TN_{1} $}. If the manifolds \textit{$N_{\theta }$} and $N_{T}$ \textit{ }(resp. $N_{\bot } $)  are  slant and invariant (resp.  anti -- invariant) submanifolds of Kenmotsu manifold $\tilde{M}$ ,then  their  warped product  semi-slant submanifolds  may given by one of the following forms:

\hspace{2.5cm} {$N_{T} \times {}_{f} N_{\theta }$, $N_{\bot } \times {}_{f} N_{\theta } $,  $N_{\theta } \times {}_{f} N_{T} $ and  $N_{\theta } \times {}_{f} N_{\bot } $.}\\
Let the warped products of type  $N_{T} \times {}_{f} N_{\theta }$ then for any $X\in TN_{T} $ and$Z\in TN_{\theta }$ we have 
\[g\left(\phi \tilde{\nabla }_{X} Z,\phi Z\right)=g\left(\tilde{\nabla }_{X} Z,Z\right).\] 
\noindent
On using equation \eqref{GrindEQ__3_2_}, we obtain 

\begin{equation} \label{GrindEQ__4_1_}
g\left(\phi \tilde{\nabla }_{X} Z,\phi Z\right)=\left(X\ln f\right)\left\| Z\right\| ^{2}.
\end{equation}
\noindent
On the other hand we have
\[\left(\tilde{\nabla }_{X} \phi \right)Z=\tilde{\nabla }_{X} \phi Z-\phi \tilde{\nabla }_{X} Z,\] 
\noindent
For any $X\in TN_{T} $ and$Z\in TN_{\theta } $. Using equation \eqref{GrindEQ__2_6_} and the fact that $\xi \in TN_{T} $, left hand side of the above equation is zero, then $\tilde{\nabla }_{X} \phi Z=\phi \tilde{\nabla }_{X} Z.$Taking the inner product with $\phi Z$ and then using equation \eqref{GrindEQ__2_10_}, we get 
\[g\left(\phi \tilde{\nabla }_{X} Z,\phi Z\right)=g\left(\tilde{\nabla }_{X} \left(TZ+NZ\right),TZ+NZ\right).\] 
\noindent
On applying \eqref{GrindEQ__2_9_}, \eqref{GrindEQ__3_2_}, Gauss and Weingarten formula we have
\[g(\phi \tilde{\nabla }_{X} Z,\phi Z)=(X\ln f)g(TZ, TZ)+g((\tilde{\nabla }_{X} N)Z, NZ)+(X\ln f)g(NZ, NZ).\] 
\noindent
On using \eqref{GrindEQ__2_16_} and \eqref{GrindEQ__2_17_}, we get

\begin{equation} \label{GrindEQ__4_2_}
g\left(\phi \tilde{\nabla }_{X} Z,\phi Z\right)=\left(X\ln f\right)\left\| Z\right\| ^{2} +g\left(\left(\tilde{\nabla }_{X} N\right)Z,NZ\right).
\end{equation}
\noindent
Using equations \eqref{GrindEQ__4_1_} and \eqref{GrindEQ__4_2_}, we get 

\begin{equation} \label{GrindEQ__4_3_}
g\left(\left(\tilde{\nabla }_{X} N\right)Z,NZ\right)=0.
\end{equation}

\noindent Since $N_{\theta } $ is proper slant submanifold of $\tilde{M}$, then equation \eqref{GrindEQ__4_3_} implies the following theorem:\textbf{}\\
\noindent \textbf{Theorem 4.1.} \textit{Let $M=N_{T} \times {}_{f} N_{\theta } $ be warped product semi -- slant submanifolds of Kenmotsu manifold $\tilde{M}$ such that $\xi \in TN_{T} $. Then $\left(\tilde{\nabla }_{X} N\right)Z$ lies in the invariant normal subbundle for all $X\in TN_{T} $ and $Z\in TN_{\theta } $ where $N_{T} $ and $N_{\theta }$ are invariant and proper slant submanifolds of $\tilde{M}$.}\\
\noindent Let the warped products of type $N_{\bot } \times {}_{f} N_{\theta } $ then for any $Z\in TN_{\bot } $ and $X\in TN_{\theta } $ we have
\[\left(\tilde{\nabla }_{X} \phi \right)Z=\tilde{\nabla }_{X} \phi Z-\phi \tilde{\nabla }_{X} Z\] 
On using equations \eqref{GrindEQ__2_7_}, \eqref{GrindEQ__2_8_}, \eqref{GrindEQ__2_10_}, \eqref{GrindEQ__2_11_}, \eqref{GrindEQ__2_6_} and the fact that $\xi \in TN_{\bot } $, we obtain\\
$-\eta (Z)TX-\eta (Z)NX=-A_{NZ} X+\nabla _{X}^{\bot } NZ-T\nabla _{X} Z-N\nabla _{X} Z-th(X, Z)-nh(X, Z).$
\noindent
Using the tangential components and using equation \eqref{GrindEQ__3_2_}, we get
\[\eta \left(Z\right)TX=A_{NZ} X+\left(Z\ln f\right)TX+th\left(X,Z\right).\] 
\noindent
Taking the product with TX in above equation and using the fact that X and TX are mutually orthogonal vector fields, then
\[\eta \left(Z\right)g\left(TX,TX\right)=g\left(A_{NZ} X,TX\right)+\left(Z\ln f\right)g\left(TX,TX\right)+g\left(th\left(X,Z\right),TX\right).\] 
\noindent
Thus from equations \eqref{GrindEQ__2_9_} and \eqref{GrindEQ__2_16_}, we get

\begin{equation} \label{GrindEQ__4_4_}
\left\{\eta \left(Z\right)-\left(Z\ln f\right)\right\}\cos ^{2} \theta \left\| X\right\| ^{2} =g\left(h\left(X,TX\right),NZ\right)-g\left(h\left(X,Z\right),NTX\right).
\end{equation}

\noindent As $N_{\theta } $ is proper slant, interchanging X by TX in above equation and taking account of equation \eqref{GrindEQ__2_15_}, we deduce that
\noindent
\begin{equation} \label{GrindEQ__4_5_} \{\eta (Z)-(Z \ln f)\}\cos^{2}\theta\| X\| ^{2}=-g(h(X, TX), NZ)+g(h(TX, Z), NX).
\end{equation}
\noindent
\noindent On adding equations \eqref{GrindEQ__4_4_} and \eqref{GrindEQ__4_5_}, we obtain\\
$\{\eta (Z)-(Z\ln f)\}\cos ^{2} \theta \| X\| ^{2} =-g(h(X, TX), NZ)+g(h(TX, Z), NX).$\\
\noindent Thus by \eqref{GrindEQ__3_13_}, the right hand side of the above equation is zero and we have\\
$\left\{\eta \left(Z\right)-\left(Z\ln f\right)\right\}\cos ^{2} \theta \left\| X\right\| ^{2} =0.$\\
\noindent This leads to the following theorem:\\
\noindent \textbf{Theorem 4.2.} \textit{Let $M=N_{\bot } \times {}_{f} N_{\theta } $ be warped product semi -- slant submanifolds of Kenmotsu manifold $\tilde{M}$ such that $\xi \in TN_{\bot }$. Then }
\begin{equation} \label{GrindEQ__4_6_}
\eta \left(Z\right)=\left(Z\ln f\right) 
\end{equation}\\
\noindent \textit{For all $Z\in TN_{\bot } $where $N_{\bot }$ and $N_{\theta } $are anti-invariant and proper slant submanifolds of $\tilde{M}$, respectively.} 

\noindent Let $M=N_{\theta } \times {}_{f} N_{T}$ be warped product semi -- slant submanifolds of a Kenmotsu manifold $\tilde{M}$ such that $\xi \in TN_{\theta } $. Then for any $X\in TN_{T} $ and $Z\in TN_{\theta }$ we have
\[\left(\tilde{\nabla }_{X} \phi \right)Z=\tilde{\nabla }_{X} \phi Z-\phi \tilde{\nabla }_{X} Z\] 
On using equations  \eqref{GrindEQ__2_7_}, \eqref{GrindEQ__2_8_}, \eqref{GrindEQ__2_10_}, \eqref{GrindEQ__2_11_}, \eqref{GrindEQ__2_6_}  and the fact that$\xi \in TN_{\bot } $, we obtain
\[-\eta \left(Z\right)\phi X-\eta \left(Z\right)\phi X=\nabla _{X} TZ+h\left(X,TZ\right)-A_{NZ} X\]
\[+\nabla _{X}^{\bot } NZ-T\nabla _{X} Z-N\nabla _{X} Z-th\left(X,Z\right)-nh\left(X,Z\right).\] 
\noindent
On comparing the tangential and normal parts we have

\begin{equation} \label{GrindEQ__4_7_}
\eta \left(Z\right)\phi X=-A_{NZ} X+\nabla _{X} TZ-T\nabla _{X} Z-th\left(X,Z\right)
\end{equation}
and

\begin{equation} \label{GrindEQ__4_8_}
\left(\tilde{\nabla }_{X} N\right)Z=nh\left(X,Z\right)-h\left(X,TZ\right).
\end{equation}
\noindent
Taking the product with NZ in \eqref{GrindEQ__4_8_} we have
\[\begin{array}{l} {g\left(\left(\tilde{\nabla }_{X} N\right)Z,NZ\right)=g\left(nh\left(X,Z\right),NZ\right)-g\left(h\left(X,TZ\right),NZ\right).} \\ {\, \, \, \, \, \, \, \, \, \, \, \, \, \, \, \, \, \, \, \, \, \, \, \, \, \, \, \, \, \, \, \, \, \, \, \, \, \, \, \, =g\left(\phi h\left(X,Z\right),\phi Z\right)-g\left(th\left(X,Z\right),TZ\right)-g\left(h\left(X,TZ\right),\phi Z\right).} \end{array}\] 
\noindent
That is,

\begin{equation} \label{GrindEQ__4_9_}
g\left(\left(\tilde{\nabla }_{X} N\right)Z,NZ\right)=-g\left(th\left(X,Z\right),TZ\right)+g\left(th\left(X,TZ\right),Z\right).
\end{equation}
\noindent
As $\theta \ne \pi /2$, then substituting Z by TZ in \eqref{GrindEQ__4_9_} and using \eqref{GrindEQ__2_15_} we obtain
\[g\left(\left(\tilde{\nabla }_{X} N\right)TZ,NTZ\right)=\cos ^{2} \theta \left\{-g\left(th\left(X,Z\right),TZ\right)+g\left(th\left(X,TZ\right),Z\right)\right\}.\] 
\noindent
 Using equation \eqref{GrindEQ__4_9_}, we get
\[g\left(\left(\tilde{\nabla }_{X} N\right)TZ,NTZ\right)=\cos ^{2} \theta \, g\left(\left(\tilde{\nabla }_{X} N\right)Z,NZ\right)\, .\, \] 
\noindent
 Hence we can state the following theorem:\\
\noindent \textbf{Theorem 4.3.} \textit{Let $M=N_{\theta } \times {}_{f} N_{T} $ be warped product semi -- slant submanifolds of Kenmotsu manifold $\tilde{M}$ such that $\xi \in TN_{\theta } $. Then }
\begin{equation} \label{GrindEQ__4_10_}
g\left(\left(\tilde{\nabla }_{X} N\right)Z,NZ\right)=\sec ^{2} \theta \, g\left(\left(\tilde{\nabla }_{X} N\right)TZ,NTZ\right)\, .\, 
\end{equation}
\noindent \textit{for all $X\in TN_{T} $ and }$Z\in TN_{\theta } $\textit{where $N_{T} $ and $N_{\theta } $are invariant and proper slant submanifolds of $\tilde{M}$,respectively.} \\
\noindent Let $M=N_{T} \times {}_{f} N_{\bot } $be warped product semi -- slant submanifolds of a Kenmotsu manifold $\tilde{M}$ such that $\xi \in TN_{\theta } $. Then for any $X\in TN_{T} $ and $Z\in TN_{\bot } $ we have
\[\left(\tilde{\nabla }_{X} \phi \right)Z=\tilde{\nabla }_{X} \phi Z-\phi \tilde{\nabla }_{X} Z\] 
\noindent
 On using equation \eqref{GrindEQ__2_6_} and the fact that $\xi \in TN_{\bot } $, we obtain
\[\tilde{\nabla }_{X} \phi Z=\phi \tilde{\nabla }_{X} Z.\] 
\noindent On using equations \eqref{GrindEQ__2_7_}, \eqref{GrindEQ__2_8_}, \eqref{GrindEQ__2_10_}, \eqref{GrindEQ__2_11_},  we get
\[-A_{NZ} X+\nabla _{X}^{\bot } NZ=T\nabla _{X} Z+N\nabla _{X} Z+th\left(X,Z\right)+nh\left(X,Z\right).\] 
\noindent
From the normal components of the above equation and formula \eqref{GrindEQ__3_2_} gives

\begin{equation} \label{GrindEQ__4_11_}
\nabla _{X}^{\bot } NZ=\left(X\ln f\right)NZ+nh\left(X,Z\right).
\end{equation}
\noindent
Taking the product in \eqref{GrindEQ__4_11_} with N$W_{1}$ for any$W_{1} \in TN_{\bot } $, we get
\[g\left(\nabla _{X}^{\bot } NZ,NW_{1} \right)=\left(X\ln f\right)g\left(NZ,NW_{1} \right)+g\left(nh\left(X,Z\right),NW_{1} \right)\] 

or, 
\[g\left(\nabla _{X}^{\bot } NZ,NW_{1} \right)=\left(X\ln f\right)g\left(\phi Z,\phi W_{1} \right)+g\left(\phi h\left(X,Z\right),\phi W_{1} \right).\] 
\noindent
Then from equation \eqref{GrindEQ__2_2_} we have

\begin{equation} \label{GrindEQ__4_12_}
g\left(\nabla _{X}^{\bot } NZ,NW_{1} \right)=\left(X\ln f\right)g\left(Z,W_{1} \right).
\end{equation}
\noindent
On the other hand for any $X\in TN_{T}$ and $Z\in TN_{\bot }$ we have
\[\left(\tilde{\nabla }_{X} N\right)Z=\nabla _{X}^{\bot } NZ-N\nabla _{X} Z.\] 
\noindent
Taking the product with N$W_{1}$, for any $W_{1} \in TN_{\bot } $ and using equation \eqref{GrindEQ__3_2_}, we get

\begin{equation} \label{GrindEQ__4_13_}
g\left(\left(\tilde{\nabla }_{X} N\right)Z,NW_{1} \right)=g\left(\nabla _{X}^{\bot } NZ,NW_{1} \right)-\left(X\ln f\right)g\left(Z,W_{1} \right).
\end{equation}
\noindent
Equations \eqref{GrindEQ__4_12_} and \eqref{GrindEQ__4_13_}, follows that
\begin{equation} \label{GrindEQ__4_14_}
g\left(\left(\tilde{\nabla }_{X} N\right)Z,NW_{1} \right)=0,
\end{equation}
\noindent
For any $X\in TN$ and $Z,W_{1} \in TN_{\bot }.$ If $W_{2} \in TN_{T} $ then using the formula \eqref{GrindEQ__2_14_}, we get
\[g\left(\left(\tilde{\nabla }_{X} N\right)Z,\phi W_{2} \right)=g\left(\nabla _{X}^{\bot } NZ,\phi W_{2} \right)-\left(X\ln f\right)g\left(Z,\phi W_{2} \right).\] 
\noindent
As $N_{T} $ is an invariant submanifold then $\phi W_{2} \in TN_{T} $for any$W_{2} \in TN_{T} $, thus using the fact that the product of tangential component with normal is zero, we obtain that 

\begin{equation} \label{GrindEQ__4_15_}
g\left(\left(\tilde{\nabla }_{X} N\right)Z,\phi W_{1} \right)=0,
\end{equation}
\noindent
For any $X,W_{2} \in TN_{T}$ and $Z\in TN_{\bot }.$
\noindent
Thus equations \eqref{GrindEQ__4_14_} and \eqref{GrindEQ__4_15_} lead to the following theorem:\\
\noindent \textbf{Theorem 4.4.}\textit{ Let $M=N_{T} \times {}_{f} N_{\bot } $ be warped product semi -- slant submanifolds of Kenmotsu manifold $\tilde{M}$ such that$\xi \in TN_{T} $. Then $\left(\tilde{\nabla }_{X} N\right)Z$ lies in the invariant normal subbundle for all $X\in TN_{T} $ and $Z\in TN_{\bot } $where $N_{T} $ and $N_{\bot }$ are invariant and anti -- invariant submanifolds of $\tilde{M}$, respectively. }

\noindent{\em Author's address:}\\
\noindent Sachin Kumar Srivastava\\
\noindent Department of Mathematics\\
\noindent Central University of Himachal Pradesh\\
\noindent Dharamshala-176215, INDIA. 
\email{\textcolor[rgb]{0.00,0.00,0.84}{sachink.ddumath@gmail.com}}

\end{document}